\documentclass [10pt]{article}
\usepackage{latexsym}   
\usepackage{amsmath,amssymb,amsthm}
\usepackage{amsbsy,eucal}

\usepackage{graphicx}
\usepackage{epstopdf}

\newtheorem{thm}{Theorem}[section]

 \newtheorem{lem}[thm]{Lemma}

 \newtheorem{defn}[thm]{Definition}

\numberwithin{equation}{section}

 \newcommand{\h}{\mathcal{H}}

\def \reals   {\mathbb{R}}

\def \z2      {\mathbb{Z}_2}
\def \Z2      {\mathbb{Z}_2}

\def \n     {\noindent}
\def \to    {\rightarrow}

\def \nin   {\hbox{\it $\in$ \kern-10.0 pt / } }
\def \nequiv   {\hbox{\it $\equiv$ \kern-11.0 pt / } }

\begin{document}

\title{A Note on Trudinger-Moser Functions and Reproducing Kernel Hilbert Spaces}

\author{\Large David G. Costa and Hossein T. Tehrani\\ \\
Department of Mathematical Sciences \\
University of Nevada Las Vegas\\
Las Vegas, NV 89154-4020\\ \\
costa@unlv.nevada.edu and tehrani@unlv.nevada.edu}

\date { }

\maketitle

\pagestyle{plain}
\thispagestyle{empty}

\vskip .2 in

\begin{abstract}
After a brief review of the definition of the Trudinger-Moser functions in dimension $N=2$ and some basic notions in the theory of ``Reproducing Kernel Hilbert Spaces (RKHS)'', we will show that there is a close connection between those two topics. More precisely, among other things, we start by considering a properly chosen multiple of the classical Trudinger-Moser family of functions in dimension $N=2$, which we denote by 
\begin{equation} 
\gamma_t (r) := \frac{1}{2\pi}\min\,\{ log \frac{1}{r}, log \frac{1}{t} \}\,,
\end{equation}
where $0 < t , r < 1$, and using the theory of RKHS we will show that 
$\gamma_t$ can be seen as  a ``bounded'' (linear) evaluation functional $u \longrightarrow u(t)$ for functions $u$ in a suitable Hilbert Space ${\cal H}$. A slightly different definition for a ''Trudinger-Moser'' type function will also be considered for $N\geq 3$.
\end{abstract}

\vskip .2in

\n {\bf Key Words:} Trudinger-Moser functions, Reproducing Kernel Hilbert Spaces.

\n {\bf AMS Subject Classification:} 46E22, 46E30, 46E35

\section{Preliminaries}  
\n The Trudinger-Moser function on the 2-dimensional unit ball in the Sobolev space $W^{1,2}_{0,r}(B_2)$ was considered in \cite{Mo} (see also \cite{Ct}) and given by the function
\begin{equation} 
\mu_s (r) := \frac{1}{\sqrt{2\pi\, log(1/s)}} \min\,\{ log\frac{1}{r}, log\frac{1}{s} \}\,,
\end{equation}
where $s, r \in (0,1)$ (see also \cite{Aot} and  \cite{At} among others). Clearly $\gamma_t$, as stated in the Abstract,  is a multiple of $\mu_s$ for $s=t$. Hence, we'll always consider the notation $\gamma_t (r)$, $0<t,r< 1$ in this note. 

The current work deals with a topic studied in \cite{Ct} involving the Trudinger-Moser function in dimension $N=2$. Although the focus in \cite{Ct} was motivated by a conjecture of Adimurthi-Struwe (see \cite{Mo}, \cite{Tr}, and (1.1) in \cite{Ct}), we will show here that the Trudinger-Moser functions can be seen through the lense of ``Reproducing Kernel Hilbert Spaces'' (see e.g. \cite{Pl} for the theory of RKHS) in all dimensions $N\geq 2$. The plan for this short note is as follows. Letting $B_N:=B_N (0,1)$ denote the unit ball centered at the of origin of the Euclidean space $\mathbb{R}^N$, we will first consider the two-dimensional case of radial functions $u\in \mathbb{W}^{1,2}_{0,r}(B_2 (0,1))$ in Section 2. 
\medskip

Section 3 is reserved for the higher dimensions $N\geq 3$, where using the fundamental solution of the Laplacian in those higher dimensions, and using an approach similar to the one we considered in the two-dimensional case, we will define the so-called ``Reproducing Kernel for $t\in (0,1)$'' in the Hilbert space $\mathbb{W}^{1,2}_{0,r}(B_N)$.
\medskip

So, let us now start by recalling the definition of a RKHS over the real numbers $\reals$, and consider in the following section (similar to (1.1)) the special multiple of the Trudinger-Moser function, namely,  
\begin{equation} 
\gamma_t (r) = \frac{1}{2\pi}\min\,\{ log \frac{1}{r}, log \frac{1}{t} \}\,,  \quad 0<t, r <1\,,
\end{equation}
where the minimum is taken for $t, r\in (0, 1)$, for each \underline{given} $t$ in $(0, 1)$. 
\smallskip

Let ${\cal F}(X,\mathbb{R})$ denote the set of all functions from $X$ to $\mathbb{R}$ with its usual operations of addition and scalar multiplication.
\smallskip

\defn One says that $\h$ is a RKHS over $\reals$ when
\smallskip

(i) $\h$ is a vector subspace of ${\cal F}(X,\mathbb{R})$,
\smallskip

(ii) $\h$ is a Hilbert space under a suitable inner-product $<\cdot, \cdot >_{\h}$ such that,
\smallskip

(iii) For every $t\in X$, the ``Evaluation Functional'' $E_t : \h \to \reals$ given by $E_t (u):= u(t)$ $\forall u\in \h$ is a bounded linear functional on $\h$. 
\medskip
\smallskip

\rm Then, it follows by the Riesz representation theorem that there exists a \underline{unique} $k_t \in \h$ such that $E_t (u):=< k_t, u >_{\h}$ for every $u\in\h$.
\medskip

\defn The vector $k_t \in \h$ is called \underline{``The Reproducing Kernel for''} $t\in X$. Hence, $E_t (u)=< k_t, u >_{\h}$ $\forall u\in\h$.

\defn The two-variable function $K : (0, 1)\times (0, 1) \to \reals$ defined by $K(t,s):= < k_t, k_s >_{\h}$ for $s, t \in (0, 1)$ is called \underline{``The Reproducing Kernel for''} 
$\h$. 
\medskip
\medskip

\section{The Two-dimensional Case}  
\n Let $B_2 := B (0,1)$ denote the unit ball centered at the origin of the 2-dimensional Euclidean space $\mathbb{R}^2$, and let $\h_2 := W^{1,2}_{0,r}(B_2)$ be the Sobolev space obtained as completion of compactly supported, radially symmetric  smooth functions  on $B_2$ under the inner-product  
\begin{equation} 
< u, v >_{\h_2} := \int_{B_2} \nabla u(x)\cdot \nabla v (x)\,dx\,.
\end{equation}
Our aim here is to show that the Sobolev space $\h_2$ (rather than the space $\mathbb{W}^{1,N}_{0,r}(B_N)$) is a RKHS if the evaluation functional $E_t$ is interpreted as follows. As we will show below, any equivalent class $[v]\in\h_2$, consisting of all functions which are equal a.e. on $B_2$, has a unique representative $v_{rad}(x)=v(|x|)$ which is a radial, locally absolutely continuous function on $(0,1]$. Therefore the evaluation functional $E_t$, defined as $E_t([v])=v_{rad}(t)$ is well defined and, as we shall see, a bounded linear functional on $\h_2$.
In Section 1, we considered the Trudinger-Moser functions introduced in \cite{Mo} arising from the original Trudinger-Moser (TM) inequality, which can be written as
\begin{equation} 
\sup\int_{B_2} e^{\alpha_2 |u|^2}\,dx < \infty\,, 
\end{equation}
where the sup is taken over the set ${\h_2}$ of all {\it radially symmetric} functions 
$u\in W^{1,2}_{0,r}(B_2)$ such that $\| \nabla u\|_2 \leq 1$, and $\alpha_2 := 4\pi$ in this 2-dimensional case (a proof of that was given by Moser in \cite{Mo}). We will show how these functions provide a family of reproducing kernels for $\h_2$. 
\medskip

To start our presentation, as in section 1, we let $X$ denote the interval $(0,1]$ and consider the space  $H_2$ consisting of \underline{locally absolutely continuous} functions $u: (0,1]\rightarrow \mathbb{R}$, vanishing at $t=1$,  such that their classical pointwise derivative $u'$ (which exists a.e. on $(0,1)$) belong to $L^2 (X)$ with the measure $d\mu= rdr$, i.e.,
\begin{equation*}
H_2:= \{ \mbox{ u is loc. abs. cont. on} \, (0,1]\ |\, u'\in L^2 (X,d\mu) \mbox{ and } u(1) = 0 \}. 
\end{equation*}
\medskip
The following lemma contains all that is needed to prove the main result of this section, Theorem 2.2\,.

\lem Let $H_2$ be the space of functions defined above. Then:
\smallskip

(i) $H_2$ is a Hilbert space under the inner-product  $< u, v >_2 = 2\pi\int_0^1 u'(r) v'(r)\,rdr$. Furthermore the linear operator $T:H_2\rightarrow \h_2$ given by $Tu (x)=u(|x|)$\, is an isometric isomorphism.                                      
\smallskip

(ii) $H_2$ is a RKHS and given $t\in X=(0,1]$, the reproduction kernel $k_t$ for the evaluation functional $E_t : H_2 \to \reals$, is given by the function $k_{t}=\gamma_t$, that is, for $u\in H_2$
\begin{equation} 
< \gamma_{t} , u >_2:=2\pi \int_0^1 \gamma_{t}^{\prime}(r) u' (r)\, rdr =u(t)\,.
\end{equation}
{\bf Proof.} Cauchy-Schwarz inequality clearly implies that  $< \cdot, \cdot >_2 $ defines an inner product on $H_2$. The remaining statements in item (i) follow if we prove that $T$ is an isometric isomorphism. 
Firstly, suppose $u\in H_2$ and let   $v:B_2\setminus \{0\}\rightarrow \mathbb{R}$ be defined by $v(x)=Tu (x)=u(|x|)$. Below, we use $(r,\theta)$ to denote the polar coordinates in $B_2$. We first show that $\nabla v$ exists and belongs to $L^2(B_2)$. In fact we prove that, for nonzero $x = (x_1,x_2) = (r\cos{\theta},r\sin{\theta})$\,,
\begin{equation} 
\nabla v (x)=u'(|x|)\frac{x}{|x|}=u'(r) (\cos{\theta},\sin{\theta})\,,
\end{equation}
from which it easily follows $||v||^2_{\h_2}=2\pi\int_0^1|u'(r)|^2rdr =||u||^2_2$. 
To prove (2.4), assume that $\phi$ is a compactly supported smooth function on $B_2$. Then:
\begin{eqnarray*}
\int_{B_2} v(x)\frac{\partial \phi}{\partial x_1}(x)\, dx & =& \int_0^{2\pi}\int_0^1 u(r)\frac{\partial \phi}{\partial x_1}(r,\theta))\, rdrd\theta \\
 & =& -\int_0^{2\pi}\int_0^1 (\int_r^{1}u'(s)\,ds) \left(\frac{\partial \phi}{\partial r}\frac{\partial r}{\partial x_1}+\frac{\partial \phi}{\partial \theta}\frac{\partial \theta}{\partial x_1}\right)\, rdrd\theta \\
 & =& -\int_0^{2\pi}\int_0^1 (\int_r^{1}u'(s)\,ds)\left(\frac{\partial \phi}{\partial r}\cos{\theta}-\frac{\partial \phi}{\partial \theta}\frac{\sin{\theta}}{r}\right)\, rdr\,d\theta \\
 & =& -\int_0^{2\pi}\int_0^1 \left(\int_0^{s}( r\frac{\partial \phi}{\partial r}\cos{\theta}-\frac{\partial \phi}{\partial \theta}\sin{\theta})\,dr\right)u'(s)\,ds\,d\theta \\
\end{eqnarray*}
where we have freely used the absolute continuity of $u$ and Fubini's theorem. On the other hand, 
$$
\int_0^{s} r\frac{\partial \phi}{\partial r}\cos{\theta}\, dr -\int_0^{s}\frac{\partial \phi}{\partial \theta}\sin{\theta}\,dr=
s\phi(s,\theta)-\frac{\partial}{\partial\theta}\left(\int_0^{s}\phi(r,\theta)\sin{\theta}\,dr\right).
$$
In addition,  
\begin{eqnarray*}
\int_0^{2\pi}\int_0^{1}u'(s)\frac{\partial}{\partial\theta}\left(\int_0^{s}\phi(r,\theta)\sin{\theta}\,dr\right)ds\,d\theta 
 & =& \int_0^{1}u'(s)\int_0^{2\pi}\frac{\partial}{\partial\theta}\left(\int_0^{s}\phi(r,\theta)\sin{\theta}\right) d\theta \, ds \\
 & =& 0\,.
\end{eqnarray*}

Therefore, we finally obtain
\begin{eqnarray*}
\int_{B_2} v(x)\frac{\partial \phi}{\partial x_1}(x)\, dx& =& \int_0^{2\pi}\int_0^{1}u'(s)\cos{\theta}\phi(r,\theta)\,sds\,d\theta \\
&=& \int_{B_2} u'(|x|)\frac{x_1}{|x|}\phi(x)dx\,.
\end{eqnarray*}
A similar calculation shows that $\frac{\partial v}{\partial x_2}(x)=u'(|x|)\frac{x_2}{|x|}$ completing the proof of (2.4). This shows that the operator $T$ maps into $\h_2$ and preserves the norm. Next,  to show that $T$ is onto, assume that $v\in\h_2$ and suppose that $(\phi_n)$ is a sequence of radially symmetric compactly supported smooth functions defined on $B_2$ that converge to $v$ in the $\h_2$ norm. Assume $x_0$ be a nonzero point in $B_2$ and let $|x_0|=r_0$. Then we can write
\begin{eqnarray*}
\phi_n(x_0)=\phi_n(r_0)=\phi_n (r_0,\theta)& = &-\int_{r_0}^{1} \phi_n'(s)\, ds\\
& = & 
-\int_{r_0}^{1} <\nabla \phi_n(s, \theta), (\cos{\theta}, \sin{\theta})> ds\\
& = & -\frac{1}{2\pi}\int_0^{2\pi} \int_{r_0}^{1} <\nabla \phi_n(s, \theta), (\cos{\theta}, \sin{\theta})> ds\, d\theta\\
& = & -\frac{1}{2\pi}\int_{|x_0|<|x|<1} <\nabla \phi_n(x), \frac{x}{|x|})>\frac{1}{|x|} dx
\end{eqnarray*}
Now taking limit as $n\rightarrow\infty$, and taking into account the fact that $\phi_n$ convergrs to $v$ in the $\h_2$ norm, it is easily seen that for $0<|x_0|$, 
$$ \frac{1}{2\pi}\int_{|x_0|<|x|<1} <\nabla \phi_n(x), \frac{x}{|x|}>\frac{1}{|x|}\, dx\rightarrow  \frac{1}{2\pi}\int_{|x_0|<|x|<1} <\nabla v(x), \frac{x}{|x|})>\frac{1}{|x|}\, dx$$
Also as  $\phi_n$ converges to $v$ in $L^p(B_2)$ for any $1\leq p<\infty$, a subsequence, still denoted by $\phi_n$, converges for a.e. $x_0$ in $B_2$ to $v(x_0)$, therefore for a.e. $x_0$ in $B_2$ we have
$$v(x_0)=\frac{1}{2\pi}\int_{|x_0|<|x|<1} <\nabla v(x), \frac{x}{|x|}>\frac{1}{|x|}\, dx,$$
Note that the right hand side of the above formula defines a function for \underline{every} $x_0$,  with $0<|x_0|\leq 1$, depending only  on $|x_0|$, i.e. for $0<r\leq 1$, the function
$$ u(r):=-\frac{1}{2\pi}\int_{r<|x|<1} <\nabla v(x), \frac{x}{|x|}>\frac{1}{|x|} \,dx,$$
 is well defined, clearly locally absolutely continuous (as $\nabla u\in L^2(B_2)$ ) and  $u(r)=v(x_0)$, a.e. on $B_2$ where $|x_0|=r$. In addition, since $u(r)=-\int_{r}^1 u'(s) ds $,  we get 
$$
u'(s)=\frac{1}{2\pi}\int_{0}^{2\pi} <\nabla v(r,\theta), (\cos{\theta}, \sin{\theta})> d\theta
$$
which, again since $\nabla v$ belongs to $L^2(B_2)$, yields $u'\in L^2(X)$ with the measure $d\mu=rdr$. On the other hand since we clearly have $T(u)=v$, the first part of the proof implies that $||v||_{\h_2}=||u||_2$. This completes the proof (i).  

As to the item (ii), let $t\in X$ be given. Recalling the definition of $\gamma_t (r)$ in (1.2), we set $k_{t} = \frac{1}{2\pi}\gamma_t$ which is clearly an element of $H_2$.
\medskip\medskip

Indeed, given an arbitrary function $u\in H_2$, we have 
\begin{eqnarray*}
< k_{t}, u >_2 &=& 2\pi \int_{0}^{1} k^{\prime}_{t}(r) u^{\prime}(r)r\,dr= \int_{0}^{1} \gamma^{\prime}_{t}(r) u^{\prime}(r) r\,dr \\
& =& \int_0^1 (log \frac{1}{r})^{\prime} \chi_{[t, 1]} u^{\prime}(r)\, r dr=-\int_0^1 \chi_{[t, 1]} u^{\prime}(r) dr \\
&=& u(t)\,,
\end{eqnarray*}
from which we conclude that the evaluation functional $E_t$ is a bounded linear functional on $H_2$ and  
\begin{equation}
\gamma_{t} \mbox{  is  ``The Reproduction Kernel for'' $t\in (0,1)$\,.} 
\end{equation}
\qed
\medskip

Considering the isomorphism given in (i) of the above lemma, we can now state the following  
\medskip

\thm Let $\h_2:=W_{0,r}^{1,2} (B_2)$ and $< u,v>_{\h_2}$ be given as in (2.1). 
Then:
\smallskip

(i) The Hilbert space $\h_2$ is a RKHS under the inner-product (2.1).
\smallskip

(ii) The reproduction kernel ${k}_{t}(x)$ for $t\in B_2\setminus \{0\}$ is given by $k_t(x)=\gamma_{|t|}(|x|)\in \h_2$.
\medskip
\section{The Higher Dimensional Case $N\geq 3$}  
\n Following our notation in section 2, let us denote by $B_N:=B(0, 1)$ the unit ball centered at the origin of the N-dimensional Euclidean space $\mathbb{\reals}^N$, and let $\h_N := W^{1,2}_{0,r}(B_N)$ denote Hilbert space of {\it radially symmetric} functions on $B_N$ under the inner-product  
\begin{equation} 
< u, v >_{\h_N} := \int_{B_N} \nabla u(x)\cdot \nabla v (x)\,dx\,.
\end{equation}
The idea here is, similar to the 2-dimensional case, to exhibit a ``reproduction kernel $k_{t}$ for $t\in B_N\setminus \{0\}$ in the space $\h_N$. As we'll show, instead of a ``logarithmic'' function in the case of $N=2$, the reproduction kernel $k_{t}$ shall involve ``power''  functions when $N\geq 3$.

We should also recall that, for any $N\geq 3$, the Trudinger-Moser inequality reads as
\begin{equation} 
\sup\int_{B_N} e^{\alpha_N |u|^{\frac{N}{N-1}}}\,dx < \infty\,, 
\end{equation}
where the sup is taken over the space of all {\it radially symmetric} functions 
$u\in W^{1,N}_{0,r}(B_N)$ such that $\| \nabla u\|_2 \leq 1$, and $\alpha_N := 
N\, {\omega_{N-1}^{1/(N-1)}}$, where $\omega_{N-1}:=|S^{N-1}|$ is the measure of the unit sphere $S^{N-1}$. 

\smallskip
As in the previous section, we let $X$ denote the interval $(0,1]$ and consider the space  $H_N$ consisting of \underline{locally absolutely continuous} functions $u: (0,1]\rightarrow \mathbb{R}$, vanishing at $t=1$,  such that their classical pointwise derivative $u'$ (which exists a.e. on $(0,1)$)  belong to $L^2 (X)$ with the measure $d\nu= r^{N-1}dr$, i.e.,
\begin{equation*}
H_N:= \{ \mbox{ u is loc. abs. cont. on} \, (0,1]\ |\, u'\in L^2 (X,d\nu) \mbox{ and } u(1) = 0 \}. 
\end{equation*}
\medskip
Similarly, as before, we can state the following
\lem 
Let $H_N$ be defined above. Then:
\smallskip

(i) $H_N$ is a RKHS under the inner-product 
$$< u, v >_N = \omega_{N-1}\int_0^1 u'(r) v'(r)\,r^{N-1}dr.$$                                            
\smallskip
\noindent

(ii) Given $t\in X=(0,1]$, the reproduction kernel $k_t$ for the evaluation functional $E_t : H_N \to \reals$, is given by the function 
\begin{equation} 
k_{t}=\gamma_{t,N}(r) : = \frac{1}{(2-N)\omega_{N-1}}\min\,\{ r^{2-N}-1, t^{2-N}-1 \}\,,\ \ 0 < r, t < 1\,
\end{equation} 
that is 
\begin{equation}< k_{t} , u >_N:=\omega_{N-1}\int_0^1 k_{t}^{\prime}(r) u' (r)\, r^{N-1}dr=u(t)\,.
\end{equation}

{\bf Proof.} As in the proof of lemma 2.1 above, one can show, using similar arguments with minor differences, that any equivalence class of functions $[v]\in\h_N$ has a unique representative $v_{rad}$ in the space $H_N$ with the correspondence providing an isomorphic isomorphism between the spaces $\h_N$ and $H_N$ with the corresponding inner products given above. Furthermore for $t\in X$, and $u\in H_N$, with $k_t$ given in (3.3), we easily calculate
\begin{eqnarray*}
< k_{t}, u >_N &=& \frac{1}{\omega_{N-1}} \int_0^1 \gamma_{t,N}^{\prime}(r) u^{\prime}(r)\,r^{N-1} dr d{\omega}\\
&= & \int_0^1 (\frac{r^{2-N}}{2-N})^{\prime}\, \chi_{[t, 1]}(r) u^{\prime}(r)\, r^{N-1} dr\\
& =& - \int_t^1 (r^{1-N})\, u^{\prime}(r)\, r^{N-1} dr\ = u(t).
\end{eqnarray*}
which completes the proof as before.

\qed
\medskip

Thus, we can state the following
\thm  
Let $\h_N := W_{0,r}^{1,2} (B_N)$, $< u, v >_{\h_N}$ be given in (3.1). 
Then: 
\medskip

(i) The Hilbert space $\h_N$ is a RKHS under the inner-product given in (3.1).
\smallskip
\noindent

(ii) The reproduction kernel ${k}_{t}(x)$ for $t\in B_N\setminus \{0\}$ is given by $k_t(x)=\gamma_{|t|,N}(|x|)\in \h_N$ where $\gamma_{t,N}$ is given in (3.3) above\,.
\medskip\medskip

We invite the reader to calculate the Reproduction Kernel for the space $\h_N$,\ i.e., the function ``$K : (0, 1)\times (0, 1) \to \reals$'' defined by 
$$
K(t,s):= < k_{t,N}, k_{s,N} >_{\h_N}, \mbox{\,for\, $s, t \in (0, 1)$\ and\ $N\geq 3$}\,.
$$

\section{Final Remarks}  
\n It is worthwhile to call attention to the following points:
\begin{itemize}
\item[1.] As already pointed out, we considered in this paper the Hilbert space $\h_2 := W_{0,r}^{1,2} (B_2)$, $N=2$, and $\h_N := W_{0,r}^{1,2} (B_N)$, $N\geq 3$ (rather than the spaces $W_{0,r}^{1,N} (B_N)$, which are only Banach spaces); 
\item[2.] We have shown that the Sobolev spaces $\h_2 := W_{0,r}^{1,2} (B_2)$, for $N=2$, and $\h_N := W_{0,r}^{1,2} (B_N)$, for $N\geq 3$, are Reproducing Kernel Hilbert Spaces (RKHSs) and determined their reproducing kernel $k_t$ for $t\in (0,1)$;
\item[3.] In [6] Moser showed that the ''sup'' in (2.2) and (3.2) were finite for $\alpha \leq \alpha_N$ (see also Trudinger in \cite{Tr} for a preliminary result);
\item[4.] In [5] Carleson-Chang proved the existence of an extremal function for the averages (which we denote below as) 
$$
Sup\, \frac{1}{|B_N|}\int_{B_N} e^{\alpha |u|^{\frac{N}{N-1}}}\,dx := C_N\,,
$$
for $N\geq 2$, and $\alpha = \alpha_N := N \omega_{N-1}^{1/(N-1)}$. Furthermore, Flucher \cite{Fl} proved the existence of an extremal function for an arbitrary bounded domain in $\mathbb{R}^2$, and Lin \cite{Li} extended the result for all $N\geq 2$;
\item[5.] To the best of our knowledge, the \underline{explicit} value of $C_N$ above (see \cite{Cc}) has not be found in general;
\item[6.] Although there are many more questions related to the existence of RKHSs of the type considered in this paper, in a subsequent paper we plan to consider the question of ''Carleson-Moser Towers'' as defined in [4] (see also references [1]-[3], among others), as well as some additional RKHSs. 
\item[7)] Motivated by (1.1) and the results in Moser [6], as well as in [2], [3] (among others), it is \underline{tempting} to call the following functions
$$
({\omega_{N-1}\,\frac{t^{2-N}-1}{N-2}})^{-1/N}\frac{1}{2-N} \min\,\{r^{2-N}-1, t^{2-N}-1 \}\,,\ \ 0 < r, t <1\,,
$$
''Moser-Trudinger Functions'' for the spaces $W^{1,2}_{0,r}(B_N)$.
\end{itemize}

\end{document}